\documentclass{amsart}
\usepackage{amssymb}
\usepackage{amsmath}

\newtheorem{theorem}{Theorem}

\newtheorem{corollary}[theorem]{Corollary}

\newtheorem{lemma}[theorem]{Lemma}

\newtheorem{proposition}[theorem]{Proposition}

\newcommand{\lm}{\operatorname{lim}}
\begin{document}

\title{Relations that preserve compact filters}
\author{Fr\'{e}d\'{e}ric Mynard}
\address{Department of Mathematical Sciences, Georgia Southern University,
Statesboro, GA 30460-8093}
\email{fmynard@georgiasouthern.edu}
\maketitle

\begin{abstract}
Many classes of maps are characterized as (possibly multi-valued) maps
preserving particular types of compact filters.
\end{abstract}

\section{Introduction}

A filter $\mathcal{F}$ on $X$ is \emph{compact} \emph{at }$A\subset X$ if
every finer ultrafilter has a limit point in $A$. As a common generalization
of compactness (in the case of a principal filter) and of convergence, it is
not surprising that the notion turned out to be very useful in a variety of
context (see for instance \cite{active}, \cite{D.comp}, \cite{cascales}
under the name of compactoid filter, \cite{pettis}, \cite{vaughan78}, \cite%
{total} under the name of total filter). The purpose of this paper is to
build on the results of \cite{active} and \cite{D.comp} to show that a large
number of classes of single and multi-valued maps classically used in
topology, analysis and optimization are instances of \emph{compact relation}%
, that is, relation that preserves compactness of filters. It is well known
(see for instance \cite{active}, \cite{cascales}) that upper semi-continuous
multivalued maps and compact valued upper semi-continuous maps are such
instances. S. Dolecki showed \cite{D.comp} that closed, countably perfect,
inversely Lindel\"{o}f and perfect maps are other examples of compact
relations. In this paper, it is shown that continuous maps as well as
various types of quotient maps (hereditarily quotient, countably biquotient,
biquotient) are also compact relations. Moreover, I\ show that maps among
these variants of quotient and of perfect maps with ranges satisfying
certain local topological properties (such as Fr\'{e}chetness, strong Fr\'{e}%
chetness and bisequentiality) can be directly characterized in similar
terms. This requires to work in the category of convergence spaces rather
than in the category of topological spaces. Therefore, I recall basic facts
on convergence spaces in the next section.

The companion paper \cite{myn.applofcompact}, which should be seen as a
sequel to the present paper, uses these characterizations to present
applications of product theorems for compact filters to theorems of
stability under product of variants of compactness, of local topological
properties (Fr\'{e}chetness and its variants, among others) and of the
classes of maps discussed above.

\section{Terminology and basic facts}

\subsection{Convergence spaces}

By a \emph{convergence space} $(X,\xi )$ I\ mean a set $X$ endowed with a
relation $\xi $ between points of $X$ and filters on $X,$ denoted $x\in
\lm_{\xi }\mathcal{F}$ or $\mathcal{F}\underset{\xi }{\rightarrow }x,$
whenever $x$ and $\mathcal{F}$ are in relation, and satisfying $\lm 
\mathcal{F}\subset \lm \mathcal{G}$ whenever $\mathcal{F}\leq \mathcal{G};$ 
$\{x\}^{\uparrow }\rightarrow x$ (\footnote{%
If $\mathcal{A}\subset 2^{X},$ $\mathcal{A}^{\uparrow }=\{B\subset X:\exists
A\in \mathcal{A},$ $A\subset B\}.$}) for every $x\in X$ and $\lm \left( 
\mathcal{F\wedge G}\right) =\lm \mathcal{F}\cap \lm \mathcal{G}$ for every
filters $\mathcal{F}$ and $\mathcal{G}$ (\footnote{%
Several different variants of these axioms have been used by various authors
under the name \emph{convergence space}.}). A map $f:(X,\xi )\rightarrow
(Y,\tau )$ is \emph{continuous }if $f(\lm_{\xi }\mathcal{F)\subset }%
\lm_{\tau }f(\mathcal{F}).$ If $\xi $ and $\tau $ are two convergences on $%
X,$ we say that $\xi $ \emph{is finer than} $\tau ,$ in symbols $\xi \geq
\tau ,$ if $Id_{X}:(X,\xi )\rightarrow (X,\tau )$ is continuous. The
category $\mathbf{Conv}$ of convergence spaces and continuous maps is
topological (\footnote{%
In other words, for every sink $\left( f_{i}:(X_{i},\xi _{i})\rightarrow
X\right) _{i\in I},$ there exists a \emph{final convergence} structure on $%
X: $ the finest convergence on $X$ making each $f_{i}$ continuous.
Equivalently, for every source $\left( f_{i}:X\rightarrow (Y_{i},\tau
_{i})\right) _{i\in I}$ there exists an\emph{\ initial convergence}: the
coarsest convergence on $X$ making each $f_{i}$ continuous.}) and
cartesian-closed (\footnote{%
In other words, for any pair $(X,\xi )$, $(Y,\tau )$ of convergence spaces,
there exists the coarsest convergence $[\xi ,\tau ]$ -called \emph{%
continuous convergence- }on the set $C(\xi ,\tau )$ of continuous functions
from $X$ to $Y$ making the evaluation map 
\begin{equation*}
ev:(X,\xi )\times \left( C(\xi ,\tau ),[\xi ,\tau ]\right) \rightarrow
(Y,\tau )
\end{equation*}%
(jointly) continuous.}).

Two families $\mathcal{A}$ and $\mathcal{B}$ of subsets of $X$ \emph{mesh, }%
in symbols $\mathcal{A}\#\mathcal{B},$ if $A\cap B\neq \emptyset $ whenever $%
A\in \mathcal{A}$ and $B\in \mathcal{B}.$ A subset $A$ of $X$ is $\xi $-%
\emph{closed} if $\lm_{\xi }\mathcal{F}\subset A$ whenever $\mathcal{F}\#A.$
The family of $\xi $-closed sets defines a topology $T\xi $ on $X$ called 
\emph{topological modification} of $\xi .$ The neighborhood filter of $x\in X
$ for this topology is denoted $\mathcal{N}_{\xi }(x)$ and the closure
operator for this topology is denoted $\mathrm{cl}_{\xi }.$ A convergence
is a topology if $x\in \lm_{\xi }\mathcal{N}_{\xi }(x).$ By definition, the
adherence of a filter (in a convergence space) is: 
\begin{equation}
\mathrm{adh}_{\xi }\mathcal{F}=\bigcup_{\mathcal{G}\#%
\mathcal{F}}\lm_{\xi }\mathcal{G}.  \label{eq:adh}
\end{equation}%
In particular, the adherence of a subset $A$ of $X$ is the adherence of its 
\emph{principal filter} $\{A\}^{\uparrow }.$ The \emph{vicinity filter} $%
\mathcal{V}_{\xi }(x)$ of $x$ for $\xi $ is the infimum of the filters
converging to $x$ for $\xi .$ A convergence $\xi $ is a \emph{pretopology}
if $x\in \lm_{\xi }\mathcal{V}_{\xi }(x).$ A convergence $\xi $ is
respectively a topology, a pretopology, a\emph{\ paratopology, }a\emph{\
pseudotopology} if $x\in \lm_{\xi }\mathcal{F}$ whenever $x\in \mathrm{adh}%
_{\xi }\mathcal{D},$ for every $\mathbb{D}$-filter $\mathcal{D}\#%
\mathcal{F}$ where $\mathbb{D}$ is respectively, \emph{the class} $\mathrm{%
cl}_{\xi }^{\natural }\left( \mathbb{F}_{1}\right) $ of principal filters of 
$\xi $-closed sets (\footnote{%
More generally, if $o:2^{X}\longrightarrow 2^{X}$ and $\mathcal{F}\subset
2^{X}$ then $o^{\natural }\mathcal{F}$ denotes $\{o(F):F\in \mathcal{F}\}$
and if $\mathbb{D}$ is a class of filters (or of family of subsets) then $%
o^{\natural }\left( \mathbb{D}\right) $ denotes $\{\mathcal{F}:\exists 
\mathcal{D}\in \mathbb{D},$ $\mathcal{F}=o^{\natural }\mathcal{D}\}.$}), 
\emph{the class} $\mathbb{F}_{1}$ of principal filters, \emph{the class }$%
\mathbb{F}_{\omega }$ of countably based filters,\emph{\ the class} $\mathbb{%
F}$ of all filters. In other words, the map $\mathrm{Adh}_{\mathbb{D}}$ 
\cite{quest2} defined by 
\begin{equation}
\lm_{\mathrm{Adh}_{\mathbb{D}}\xi }\mathcal{F}=\bigcap_{%
\mathbb{D}\backepsilon \mathcal{D}\#\mathcal{F}}\mathrm{adh}_{\xi }%
\mathcal{D}  \label{eq:Adh}
\end{equation}%
is the (restriction to objects of the) reflector from \textbf{Conv} onto the
full subcategory of respectively topological, pretopological,
paratopological and pseudotopological spaces when $\mathbb{D}$ is
respectively, the class\emph{\ }$\mathrm{cl}_{\xi }^{\natural }\left( 
\mathbb{F}_{1}\right) $, $\mathbb{F}_{1}$, $\mathbb{F}_{\omega }$ and $%
\mathbb{F}$.

A convergence space is \emph{first-countable }if whenever $x\in \lm 
\mathcal{F},$ then there exists a countably based filter $\mathcal{H}\leq 
\mathcal{F}$ such that $x\in \lm \mathcal{H}.$ Of course, a topological
space is first-countable in the usual sense if and only if it is
first-countable as a convergence space. Analogously, a convergence space is
called \emph{sequentially based }if whenever $x\in \lm \mathcal{F},$ there
exists a sequence $(x_{n})_{n\in \omega }\leq \mathcal{F}$ (\footnote{%
From the viewpoint of convergence, there is no reason to distinguish between
a sequence and the filter generated by the family of its tails. Therefore,
in this paper, sequences are identified to their associated filter and I
will freely treat sequences as filters. Hence the notation $(x_{n})_{n\in
\omega }\leq \mathcal{F}.$}) such that $x\in \lm (x_{n})_{n\in \omega }.$

A class of filters $\mathbb{D}$ (under mild conditions on $\mathbb{D}$)
defines a reflective subcategory of \textbf{Conv} (and the associated
reflector) via $\left( \ref{eq:Adh}\right) $. Dually, it also defines (under
mild conditions on $\mathbb{D}$) the coreflective subcategory of \textbf{Conv%
} of $\mathbb{D}$-\emph{based convergence spaces }\cite{quest2}, and the
associated (restriction to objects of the) coreflector $\mathrm{Base}_{%
\mathbb{D}}$ is 
\begin{equation}
\lm_{\mathrm{Base}_{\mathbb{D}}\xi }\mathcal{F}=\bigcup_{%
\mathbb{D}\backepsilon \mathcal{D}\leq \mathcal{F}}\lm_{\xi }%
\mathcal{D}.  \label{eq:Base}
\end{equation}

For instance, if $\mathbb{D}=\mathbb{F}_{\omega }$ is the class of countably
based filter, then $\mathrm{Base}_{\mathbb{D}}$ is the coreflector on
first-countable convergence spaces. If\emph{\ }$\mathbb{D}$ is the class $%
\mathbb{E}$ of filters generated by sequences, then $\mathrm{Base}_{\mathbb{%
D}}$ is the coreflector on sequentially based convergences.

\subsection{Local properties and special classes of filters}

Recall that a topological space is \emph{Fr\'{e}chet }(respectively, \emph{%
strongly Fr\'{e}chet}) if whenever $x$ is in the closure of a subset $A$
(respectively, $x$ is in the intersection of closures of elements of a
decreasing sequence $(A_{n})_{n}$ of subsets of $X)$ there exists a sequence 
$(x_{n})_{n\in \omega }$ of elements of $A$ (respectively, such that $%
x_{n}\in A_{n}$) such that $x\in \lm (x_{n})_{n\in \omega }.$ In other
words, if $x$ is in the adherence of a principal (resp. countably based)
filter, then there exists a sequence meshing with that filter that converges
to $x.$ These are special cases of the following general notion, defined for
convergence spaces.

Let $\mathbb{D}$ and $\mathbb{J}$ be two classes of filters. A convergence
space $(X,\xi )$ is called $(\mathbb{J}/\mathbb{D})$-\emph{accessible} if 
\begin{equation*}
\mathrm{adh}_{\xi }\mathcal{J}\subset \mathrm{adh}_{%
\mathrm{Base}_{\mathbb{D}}\xi }\mathcal{J},
\end{equation*}%
for every $\mathcal{J}\in \mathbb{J}.$ When $\mathbb{D=F}_{\omega }$ and $%
\mathbb{J}$ is respectively the class $\mathbb{F},$ $\mathbb{F}_{\omega }$
and $\mathbb{F}_{1},$ then $(\mathbb{J}/\mathbb{D})$-accessible topological
spaces are respectively bisequential, strongly Fr\'{e}chet and Fr\'{e}chet
spaces. Analogously, if $\mathbb{D}$ is the class of filters generated by
long sequences (of arbitrary length) and $\mathbb{J}=\mathbb{F}_{1}$ then $(%
\mathbb{J}/\mathbb{D})$-accessible topological spaces are radial spaces. We
use the same names for these instances of $(\mathbb{J}/\mathbb{D})$%
-accessible \emph{convergence spaces }(see \cite{quest2} for details).

A filter $\mathcal{F}$ is called $\mathbb{J}$ to\emph{\ }$\mathbb{D}$\emph{\
meshable-refinable}, in symbol $\mathcal{F}\in (\mathbb{J}/\mathbb{D}%
)_{\#\geq }$, if 
\begin{equation*}
\mathcal{J}\in \mathbb{J},\text{ }\mathcal{J}\#\mathcal{F}\Longrightarrow
\exists \mathcal{D}\in \mathbb{D},\text{ }\mathcal{D}\#\mathcal{J}\text{ and 
}\mathcal{D}\geq \mathcal{F}.
\end{equation*}%
It follows immediately from the definitions that a\emph{\ topological} space
is $(\mathbb{J}/\mathbb{D})$-accessible if and only if every neighborhood
filter is $\mathbb{J}$ to $\mathbb{D}$ meshable-refinable, and more
generally that:

\begin{theorem}
\label{th:accessible}Let $\mathbb{D}$ and $\mathbb{J}$ be two classes of
filters.

\begin{enumerate}
\item A convergence space $(X,\xi )$ is $(\mathbb{J}/\mathbb{D})$-accessible
if and only if $\xi \geq \mathrm{Adh}_{\mathbb{J}}\mathrm{Base}_{\mathbb{D}%
}\xi $;

\item If $\xi =\mathrm{Base}_{(\mathbb{J}/\mathbb{D})_{\#\geq }}\xi ,$ then 
$\xi $ is $(\mathbb{J}/\mathbb{D})$-accessible. If moreover $\xi $ is
pretopological (in particular topological) then the converse is true.
\end{enumerate}
\end{theorem}

The following gathers the most common cases of $(\mathbb{J}/\mathbb{D})$%
-accessible (topological) spaces and $(\mathbb{J}/\mathbb{D})_{\#\geq }$%
-filters when $\mathbb{D=F}_{\omega }.$ Denote by $\mathbb{F}_{\wedge \omega
}$ the class of \emph{countably deep} (\footnote{%
A filter $\mathcal{F}$ is \emph{countably deep }if $\bigcap \mathcal{A}\in 
\mathcal{F}$ whenever $\mathcal{A}$ is a countable subfamily of $\mathcal{F}%
. $}) filters. The names for $(\mathbb{J}/\mathbb{F}_{\omega })_{\#\geq }$%
-filters come from the fact that a topological space is $(\mathbb{J}/\mathbb{%
F}_{\omega })$-accessible if and only if every neighborhood filter is a $(%
\mathbb{J}/\mathbb{F}_{\omega })_{\#\geq }$-filter.\bigskip

\begin{center}
\begin{tabular}{|c|c|c|}
\hline
class $\mathbb{J}$ & $(\mathbb{J}/\mathbb{F}_{\omega })$-accessible space & $%
(\mathbb{J}/\mathbb{F}_{\omega })_{\#\geq }$-filter \\ \hline
$\mathbb{F}$ & bisequential \cite{quest} & bisequential \\ \hline
$\mathbb{F}_{\omega }$ & strongly Fr\'{e}chet or countably bisequential \cite%
{quest} & strongly Fr\'{e}chet \\ \hline
$(\mathbb{F}_{\omega }/\mathbb{F}_{\omega })_{\#\geq }$ & productively Fr%
\'{e}chet\cite{JM} & productively Fr\'{e}chet \\ \hline
$\mathbb{F}_{\wedge \omega }$ & weakly bisequential \cite{arh.biseq} & 
weakly bisequential \\ \hline
$\mathbb{F}_{1}$ & Fr\'{e}chet \cite{quest} & Fr\'{e}chet \\ \hline
\end{tabular}

Table 1
\end{center}

\subsection{Compactness}

Let $\mathbb{D}$ be a class of filters on a convergence space $(X,\xi )$ and
let $\mathcal{A}$ be a family of subsets of $X$. A filter $\mathcal{F}$ is $%
\mathbb{D}$-\emph{compact }at $\mathcal{A}$ (for $\xi )$ \cite{comp} if (%
\footnote{%
Notice that (\ref{eq:Dcompactoid}) makes sense not only for a filter but for
a general family $\mathcal{F}$ of subsets of $X$. Such general compactoid
families play an important role for instance in \cite{DGL.kur}.}) 
\begin{equation}
\mathcal{D\in }\mathbb{D},\mathcal{D}\#\mathcal{F}\Longrightarrow \mathrm{%
adh}_{\xi }\mathcal{D}\#\mathcal{A}.  \label{eq:Dcompactoid}
\end{equation}

Notice that a subset $K$ of a convergence space $X$ (in particular of a
topological space) is respectively compact, countably compact, Lindel\"{o}f
if $\{K\}^{\uparrow }$ is $\mathbb{D}$-compact at $\{K\}$ if $\mathbb{D}$ is
respectively, the class $\mathbb{F}$ of all, $\mathbb{F}_{\omega }$ of
countably based, $\mathbb{F}_{\wedge \omega }$ of countably deep filters. On
the other hand

\begin{theorem}
\label{pro:AdhD} Let $\mathbb{D}$ be a class of filters. A filter $\mathcal{F%
}$ is $\mathbb{D}$-compact at $\{x\}$ for $\xi $ if and only if 
\begin{equation*}
x\in \lm_{\mathrm{Adh}_{\mathbb{D}}\xi }\mathcal{F}.
\end{equation*}
\end{theorem}

In particular, if $\xi $ is a topology, then $x\in \lm \mathcal{F}$ if and
only if $\mathcal{F}$ is compact at $\{x\}$ if and only if $\mathcal{F}$ is $%
\mathbb{F}_{1}$-compact at $\{x\}.$

For a topological space $X,$ a subset $K$ is compact if and only if every
open cover of $K$ has a finite subcover of $K,$ if and only if every filter
on $K$ has adherent points in $K.$ In contrast, for general convergence
spaces, the definition of compactness in terms of covers (cover-compactness)
and in terms of filters (compactness) are different. If $(X,\xi )$ is a
convergence space, a family $\mathcal{S}\subset 2^{X}$ is a \emph{cover of }$%
K\subset X$ if every filter converging to a point of $K$ contains an element
of $\mathcal{S}.$ Hence a subset $K$ of a convergence space is called \emph{%
cover-(countably ) compact }if every (countable) cover of $K$ has a finite
subcover. It is easy to see that a cover-compact convergence is compact, but
in general not conversely. For instance, in a pseudotopological but not
pretopological convergence, points are compact, but not cover-compact. 

Notice that in this definition, we can assume the original cover $\mathcal{S}
$ to be stable under finite union, in which case we call $\mathcal{S}$ an 
\emph{additive cover}. The family $\mathcal{S}_{c}$ of complements of
elements of an additive cover $\mathcal{S}$ is a filter-base on $X$ with
empty adherence. Hence $K$ is cover-(countably) compact if every (countable)
additive cover of $K$ has an element that is a cover of $K,\ $or
equivalently, if every (countable) filter-base with no adherence point in $K$
has an element with no adherence point in $K.$ In other words, $K$ is
cover-(countably) compact if every (countably based) filter whose every
member has adherent points in $K$, has adherent points in $K.$

More generally, we will need the following characterization of
cover-compactness in terms of filters \cite{D.comp}. Let $\mathbb{D}$ and $%
\mathbb{J}$ be two classes of filters. A filter $\mathcal{F}$ is $(\mathbb{D}%
/\mathbb{J})$-\emph{compact at} $\mathcal{B}$ if 
\begin{equation*}
\mathcal{D}\in \mathbb{D},\text{ }\forall \mathcal{J}\in \mathbb{J},\text{ }%
\mathcal{J}\leq \mathcal{D},\text{ }\mathrm{adh}\mathcal{J}\#\mathcal{%
F\Longrightarrow }\mathrm{adh}\mathcal{D}\#\mathcal{B}.
\end{equation*}

It is clear than if $\mathcal{F}$ is $(\mathbb{D}/\mathbb{F}_{1})$-compact
(at $\mathcal{B)},$ then it is $\mathbb{D}$-compact (at $\mathcal{B)}.$ More
precisely, we have the following relationship between $(\mathbb{D}/\mathbb{F}%
_{1})$-compactness and $\mathbb{D}$-compactness (which could be deduced from
the results of \cite[section 8]{D.comp})

\begin{proposition}
Let $\mathbb{D}$ be a class of filters on a convergence space $(X,\xi ).$ A
filter $\mathcal{F}$ is $(\mathbb{D}/\mathbb{F}_{1})$-compact at $\mathcal{B}
$ if and only if $\mathcal{V}_{\xi }(\mathcal{F})=\bigcup_{F\in 
\mathcal{F}}\bigcap_{x\in F}\mathcal{V}_{\xi }(x)$ is $\mathbb{D}$%
-compact at $\mathcal{B}.$
\end{proposition}

\begin{proof}
By definition, $\mathcal{F}$ is $(\mathbb{D}/\mathbb{F}_{1})$-compact at $%
\mathcal{B}$ if and only if 
\begin{equation*}
\mathcal{D}\in \mathbb{D}:\left( \mathrm{adh}_{\xi }^{\natural }%
\mathcal{D}\right) \#\mathcal{F\Longrightarrow }\mathrm{adh}\mathcal{D}\#%
\mathcal{B}.
\end{equation*}%
It is easy to verify that $\left( \mathrm{adh}_{\xi }^{\natural }%
\mathcal{D}\right) \#\mathcal{F}$ if and only if $\mathcal{D}\#\mathcal{V}%
_{\xi }(\mathcal{F}),$ which concludes the proof.
\end{proof}

Calling a convergence $\xi $ \emph{pretopologically diagonal}, or $P$-\emph{%
diagonal}, if $\lm_{\xi }\mathcal{F}\subset \lm_{\xi }\mathcal{V}_{\xi }(%
\mathcal{F})$ for every filter $\mathcal{F},$ we obtain the following
result, which is a particular case of a combination of Propositions 8.1 and
8.3 and of Theorem 8.2 in \cite{D.comp}, even though the assumption that $%
\mathrm{adh}_{\xi }^{\natural }\mathbb{D\subset D}$ seems to be
erroneously missing in \cite{D.comp}.

\begin{corollary}
\label{cor:Pdiag} If $\xi $ is $P$-diagonal (in particular if $\xi $ is a
topology) and if $\mathrm{adh}_{\xi }^{\natural }\mathbb{D\subset D%
}$, then $(\mathbb{D}/\mathbb{F}_{1})$-compactness amounts to $\mathbb{D}$%
-compactness for $\xi .$
\end{corollary}

\begin{proof}
Assume that $\mathcal{F}$ is $\mathbb{D}$-compact(at $\mathcal{B}$). To show
that it is $(\mathbb{D}/\mathbb{F}_{1})$-compact(at $\mathcal{B}$), we only
need to show that $\mathcal{V}_{\xi }(\mathcal{F})$ is $\mathbb{D}$%
-compact(at $\mathcal{B}$). But $\mathcal{D}\#\mathcal{V}_{\xi }(\mathcal{F}%
) $ if and only if $\left( \mathrm{adh}_{\xi }^{\natural }\mathcal{%
D}\right) \#\mathcal{F}.$ Therefore, $\mathrm{adh}_{\xi }(\mathrm{%
adh}_{\xi }^{\natural }\mathcal{D})\#\mathcal{B}$ because $\mathrm{%
adh}_{\xi }^{\natural }\mathcal{D}\in \mathbb{D}$. Now, $x\in 
\mathrm{adh}_{\xi }(\mathrm{adh}_{\xi }^{\natural }%
\mathcal{D})$ if there exists a filter $\mathcal{G}\#\mathrm{adh}%
_{\xi }^{\natural }\mathcal{D}$ with $x\in \lm_{\xi }\mathcal{G}.$
Note that $\mathcal{V}_{\xi }(\mathcal{G})\#\mathcal{D}$ and that, by $P$%
-diagonality, $x\in \lm_{\xi }\mathcal{V}_{\xi }(\mathcal{G}).$ Hence $%
\mathrm{adh}_{\xi }(\mathrm{adh}_{\xi }^{\natural }%
\mathcal{D})\subset \mathrm{adh}_{\xi }\mathcal{D}$ and $\mathrm{%
adh}_{\xi }\mathcal{D}\#\mathcal{B}.$
\end{proof}

In some sense, the converse is true:

\begin{proposition}
If $\xi =\mathrm{Adh}_{\mathbb{D}}\xi $ and $\mathbb{D}$%
-compactness implies $(\mathbb{D}/\mathbb{F}_{1})$-compactness in $\xi ,$
then $\xi $ is $P$-diagonal.
\end{proposition}

\begin{proof}
If $x\in \lm_{\xi }\mathcal{F}$ then $\mathcal{F}$ is $\mathbb{D}$-compact
at $\{x\},$ hence $(\mathbb{D}/\mathbb{F}_{1})$-compact at $\{x\}.$ Since $%
\xi =\mathrm{Adh}_{\mathbb{D}}\xi ,$ we only need to show that $%
x\in \mathrm{adh}_{\xi }\mathcal{D}$ whenever $\mathcal{D}$ is a $\mathbb{D}
$-filter meshing with $\mathcal{V}_{\xi }(\mathcal{F}).$ For any such $%
\mathcal{D},$ we have $\mathrm{adh}_{\xi }^{\natural }\mathcal{D}\#%
\mathcal{F}$ so that $x\in \mathrm{adh}_{\xi }\mathcal{D}$ because $%
\mathcal{F}$ is $(\mathbb{D}/\mathbb{F}_{1})$-compact at $\{x\}.$
\end{proof}

\subsection{Contour filters}

If $\mathcal{F}$ is a filter on $X$ and $\mathcal{G}:X\rightarrow \mathbb{F}%
X $ then the \emph{contour of} $\mathcal{G}$ \emph{along }$\mathcal{F}$ is
the filter on $X$ defined by%
\begin{equation*}
\int_{\mathcal{F}}\mathcal{G}=\bigvee_{F\in \mathcal{F}%
}\bigwedge_{x\in F}\mathcal{G}(x).
\end{equation*}

This type of filters have been used in many situations, among others by Frol%
\'{\i}k under the name of sum of filters for a ZFC proof of the
non-homogeneity of the remainder of $\beta \mathbb{N}$ \cite{Frolik}, by C.
H. Cook and H. R. Fisher \cite{cook} under the name of compression operator
of $\mathcal{F}$ relative to $\mathcal{G},$ by H. J. Kowalsky \cite{kowalsky}
under the name of diagonal filter, and after them by many other authors to
characterize topologicity and regularity of convergence spaces. To
generalize this construction, I need to reproduce basic facts on cascades
and multifilters. Detailed information on this topic can be found in \cite%
{cascades}.

If $(W,\sqsubseteq )$ is an ordered set, then we write 
\begin{equation*}
W(w)=\{x\in W:w\sqsubseteq x\}.
\end{equation*}%
An ordered set $(W,\sqsubseteq )$ is \textit{well-capped} if its every non
empty subset has a maximal point (\footnote{%
In other words, a well-capped ordered set is a well-founded ordered set for
the inverse order.}). Each well-capped set admits the (upper) \textit{rank}
to the effect that $r(w)=0$ if $w\in \max W$, and for $r(w)>0$, 
\begin{equation*}
r(w)=r_{W}(w)=\sup_{v\sqsupset w}(r(v)+1).
\end{equation*}

A well-capped tree with least element is called a \textit{cascade}; the
least element of a cascade $V$ is denoted by $\varnothing =\varnothing _{V}$
and is called the \textit{estuary} of $V$. It follows from the definition
that each element of a cascade is of finite length. The \textit{rank} of a
cascade is by definition the rank of its estuary. A cascade is a \textit{%
filter cascade} if its every (non maximal) element is a filter on the set of
its immediate successors.

A map $\Phi :V\setminus \{\varnothing _{V}\}\rightarrow X$, where $V$ is a
cascade, is called a \textit{multifilter on} $X$. We talk about a
multifilter $\Phi :V\rightarrow X$ under the understanding that $\Phi $ is
not defined at $\varnothing _{V}$.

A couple $(V,\Phi _{0})$ where $V$ is a cascade and $\Phi _{0}:\max
V\rightarrow A$ is a called a \textit{perifilter on} $A$. In the sequel we
will consider $V$ implicitly talking about a perifilter $\Phi _{0}$. If $%
\Phi |_{\max V}=\Phi _{0}$, then we say that the multifilter $\Phi $ is an 
\textit{extension} of the perifilter $\Phi _{0}$. The rank of a multifilter
(perifilter) is, by definition, the rank of the corresponding cascade. If${\ 
}\mathbb{D}$ is a class of filters, we call ${\ }\mathbb{D}$-\emph{%
multifilter} a multifilter with a cascade of ${\ }\mathbb{D}$-filters as
domain.

The contour of a multifilter $\Phi :V\rightarrow X$ depends entirely on the
underlying cascade $V$ and on the restriction of $\Phi $ to $\max V$, hence
on the corresponding perifilter $(V,\Phi |_{\max V})$. Therefore we shall
not distinguish between the contours of multifilters and of the
corresponding perifilters. The \textit{contour} of $\Phi :W\rightarrow X$ is
defined by induction to the effect that $\int \Phi =\Phi _{\natural
}(\varnothing _{W})$ if $r(\Phi )=1$, and (\footnote{$\Phi (v)$ is the image
by $\Phi $ of $v$ treated as a point of $V$, while $\Phi _{\natural }(v)$ is
the filter generated by $\{\Phi (F):F\in v\}$.})%
\begin{equation*}
\int \Phi =\int_{\varnothing _{W}}\left( \int \Phi |_{W(.)}\right)
\end{equation*}%
otherwise. With each class $\mathbb{D}$ of filters we associate the class $%
\int \mathbb{D}$ of all $\mathbb{D}$-contour filters, i.~e., the contours of 
$\mathbb{D}$-multifilter.

\begin{lemma}
\label{lem:contourcompose}Let $\mathbb{D}$ and $\mathbb{J}$ be two classes
of filters. If $\mathbb{D}$ is a $\mathbb{J}$-composable class of filters,
then $\int \mathbb{D}$ is also $\mathbb{J}$-composable.
\end{lemma}

\begin{proof}
We proceed by induction on the rank of a $\mathbb{D}$-multifilter. The case
of rank 1 is simply $\mathbb{J}$-composability of $\mathbb{D}.$ Assume that
for each $\mathbb{D}$-multifilter $\Phi $ on $X$ of rank $\beta $ smaller
than $\alpha $ and each $\mathbb{J}$-filter $\mathcal{J}$ on $X\times Y,$
the filter $\mathcal{J}(\int \Phi )$ is the contour of some $\mathbb{D}$%
-multifilter on $Y.$ Consider now a $\mathbb{D}$-multifilter $\left( \Phi
,V\right) $ on $X$ of rank $\alpha $ and a $\mathbb{J}$-filter $\mathcal{J}$
on $X\times Y.$ Then%
\begin{equation*}
\int \Phi =\int_{\varnothing _{V}}\left( \int \Phi |_{V(.)}\right)
=\bigvee_{F\in \varnothing _{V}}\bigwedge_{v\in F}\int \Phi
|_{V(v)},
\end{equation*}%
and%
\begin{equation*}
\mathcal{J}\left( \int \Phi \right) =\bigvee_{F\in \varnothing
_{V}}\bigwedge_{v\in F}\mathcal{J}\left( \int \Phi |_{V(v)}\right) .
\end{equation*}%
As each $\Phi |_{V(v)}$ is a multifilter of rank smaller than $\alpha ,$
each $\mathcal{J}\left( \int \Phi |_{V(v)}\right) $ is a $\left( \int 
\mathbb{D}\right) $-filter. Moreover $\varnothing _{V}$ is a $\mathbb{D}$%
-filter, so that $\mathcal{J}\left( \int \Phi \right) $ is a contour of $%
\left( \int \mathbb{D}\right) $-filters along a $\mathbb{D}$-filter, hence a 
$\left( \int \mathbb{D}\right) $-filter.
\end{proof}

\section{Compact relations}

A \emph{relation }$R:(X,\xi )\rightrightarrows (Y,\tau )$\emph{\ is }$%
\mathbb{D}$\emph{-compact }if for every subset $A$ of $X$ and every filter $%
\mathcal{F}$ that is $\mathbb{D}$-compact at $A,$ the filter $R\mathcal{F}$
is $\mathbb{D}$-compact at $RA$.

If $\mathbb{D}$ and $\mathbb{J}$ are two classes of filters, we say that $%
\mathbb{J}$ is $\mathbb{D}$-\emph{composable} if for every $X$ and $Y,$ the
(possibly degenerate) filter $\mathcal{HF}$=$\{HF:$ $H\in \mathcal{H},F\in 
\mathcal{F}\}^{\uparrow }$ (\footnote{$HF=\{y\in Y:(x,y)\in H$ and $x\in
F\}. $}) belongs to $\mathbb{J}(Y)$ whenever $\mathcal{F\in }\mathbb{J}(X)$
and $\mathcal{H}\in \mathbb{D}(X\times Y),$ with the convention that every
class of filters contains the degenerate filter. If a class $\mathbb{D}$ is $%
\mathbb{D}$-composable, we simply say that $\mathbb{D}$ is \emph{composable}%
. Notice that 
\begin{equation}
\mathcal{H}\#\left( \mathcal{F\times G}\right) \Longleftrightarrow \mathcal{%
HF}\#\mathcal{G}\Longleftrightarrow \mathcal{H}^{-}\mathcal{G}\#\mathcal{F},
\label{eq:grill}
\end{equation}%
where $\mathcal{H}^{-}\mathcal{G}$=$\{H^{-}G=\{x\in X:(x,y)\in H$ and $y\in
G\}:H\in \mathcal{H},G\in \mathcal{G}\}^{\uparrow }.$

\begin{proposition}
\label{pro:pointsuffice} If $\mathbb{D}$ is $\mathbb{F}_{1}$-composable,
then $R:(X,\xi )\rightrightarrows (Y,\tau )$ is $\mathbb{D}$-compact if and
only if $R\mathcal{F}$ is $\mathbb{D}$-compact at $Rx$ whenever $x\in
\lm_{\xi }\mathcal{F}.$
\end{proposition}

\begin{proof}
Only the \textquotedblright if\textquotedblright\ part needs a proof, so
assume that $R\mathcal{F}$ is $\mathbb{D}$-compact at $Rx$ whenever $x\in
\lm_{\xi }\mathcal{F},$ and consider a filter $\mathcal{G}$ on $X$ which is 
$\mathbb{D}$-compact at $A.$ Let $\mathcal{D}\#R\mathcal{G}$ be a $\mathbb{D}
$-filter on $Y$. Then $R^{-}\mathcal{D}\#\mathcal{G}$ so that there exists $%
x\in A\cap \mathrm{adh}_{\xi }R^{-}\mathcal{D}.$ Therefore, there exists $%
\mathcal{U}\#R^{-}\mathcal{D}$ such that $x\in \lm_{\xi }\mathcal{U}.$ By
assumption, $R\mathcal{U}$ is $\mathbb{D}$-compact at $Rx\subset RA.$ Since $%
\mathcal{D}\#R\mathcal{U},$ the filter $\mathcal{D}$ has adherent points in $%
Rx$ hence in $RA.$
\end{proof}

\begin{corollary}
\label{cor:continuous} Let $\mathbb{D}$ be an $\mathbb{F}_{1}$-composable
class of filters and let $f:(X,\xi )\rightarrow (Y,\tau )$ with $\tau =%
\mathrm{Adh}_{\mathbb{D}}\tau $. The following are equivalent:

\begin{enumerate}
\item $f$ is continuous;

\item $f$ is a compact relation;

\item $f$ is a $\mathbb{D}$-compact relation.
\end{enumerate}
\end{corollary}

\begin{proof}
$(1\Longrightarrow 2).$If $x\in \lm_{\xi }\mathcal{F},$ then $f(x)\in
\lm_{\tau }f(\mathcal{F})$ so that $f(\mathcal{F})$ is compact at $f(x)$
and $f$ is a compact relation by Proposition \ref{pro:pointsuffice}. $%
(2\Longrightarrow 3)$ is obvious and $(3\Longrightarrow 1)$ follows from
Proposition \ref{pro:AdhD}.
\end{proof}

In particular, $\mathbb{F}_{1}$-compact (equivalently compact) maps between
pretopological spaces (in particular between topological spaces) are exactly
the continuous ones.

Notice that when $\mathbb{D}$ contains the class of principal filters, then
a $\mathbb{D}$-compact relation $R$ is $\mathbb{F}_{1}$-compact and $Rx$ is $%
\mathbb{D}$-compact for each $x$ in the domain of $R,$ because $%
\{x\}^{\uparrow }$ is $\mathbb{D}$-compact at $\{x\}.$ When the cover and
filter versions of compactness coincide (in particular, in a topological
space), the converse is true:

\begin{proposition}
Let $\mathbb{D}$ be an $\mathbb{F}_{1}$-composable class of filters. If $%
R:(X,\xi )\rightrightarrows (Y,\tau )$ is an $\mathbb{F}_{1}$-compact
relation and if $Rx$ is $(\mathbb{D}/\mathbb{F}_{1})$-compact in $\tau $ for
every $x\in X,$ then $R$ is $\mathbb{D}$-compact.
\end{proposition}

\begin{proof}
Using Proposition \ref{pro:pointsuffice}, we need to show that $R\mathcal{F}$
is $\mathbb{D}$-compact at $Rx$ whenever $x\in \lm_{\xi }\mathcal{F}.$
Consider a $\mathbb{D}$-filter $\mathcal{D}\#R\mathcal{F}.$ Then, $\mathrm{%
adh}_{\tau }D\#Rx$ for every $D\in \mathcal{D}$ so that $\mathrm{adh}_{\tau
}\mathcal{D}\#Rx,$ because $Rx$ is $\frac{\mathbb{D}}{\mathbb{F}_{1}}$%
-compact.
\end{proof}

In view of Corollary \ref{cor:Pdiag}, we obtain:

\begin{corollary}
\label{cor:eq} Let $\mathbb{D}$ be an $\mathbb{F}_{1}$-composable class of
filters such that $\mathrm{adh}_{\tau }^{\natural }\mathbb{D(\tau
)\subset D(\tau )}$ and let $\tau $ be a $P$-diagonal convergence (for
instance a topology). Then $R:(X,\xi )\rightrightarrows (Y,\tau )$ is $%
\mathbb{D}$-compact if and only if it is $\mathbb{F}_{1}$-compact and $Rx$
is $\mathbb{D}$-compact in $\tau $ for every $x\in X.$
\end{corollary}

An immediate corollary of \cite[Theorem 8.1]{cascades} is that for a
topology, $\mathbb{D}$-compactness amounts to $\left( \int \mathbb{D}\right) 
$-compactness, provided that $\mathbb{D}$ is a composable class of filters.
However, the proof of \cite[Theorem 8.1]{cascades} only uses $\mathbb{F}_{1}$%
-composability of $\mathbb{D}$. Consequently,

\begin{corollary}
\label{cor:fibers} Let $\mathbb{D}$ be an $\mathbb{F}_{1}$-composable class
of filters and let $\tau $ be a topology such that $\mathrm{adh}%
_{\tau }^{\natural }\mathbb{D(\tau )\subset D(\tau )}$. Let $%
R:(X,\xi )\rightrightarrows (Y,\tau )$ be a relation. The following are
equivalent:

\begin{enumerate}
\item $R$ is $\mathbb{D}$-compact;

\item $R$ is $\mathbb{F}_{1}$-compact and $Rx$ is $\mathbb{D}$-compact in $%
\tau $ for every $x\in X;$

\item $R$ is $\mathbb{F}_{1}$-compact and $Rx$ is $\left( \int \mathbb{D}%
\right) $-compact in $\tau $ for every $x\in X;$

\item $R$ is $\left( \int \mathbb{D}\right) $-compact.
\end{enumerate}
\end{corollary}

\begin{proof}
$\left( 1\Longleftrightarrow 2\right) $ and $\left( 3\Longleftrightarrow
4\right) $ follow from Corollary \ref{cor:eq} and $\left(
1\Longleftrightarrow 4\right) $ follows from \cite[Theorem 8.1]{cascades}.
\end{proof}

The observation that perfect, countably perfect and closed maps can be
characterized as $\mathbb{D}$-compact relations is due to S. Dolecki \cite[%
section 10]{D.comp}. Recall that a surjection $f:X\rightarrow Y$ between two
topological spaces is \emph{closed} if the image of a closed set is closed
and\emph{\ perfect }(resp. \emph{countably perfect}, resp. \emph{inversely
Lindel\"{o}f}) if it is closed with compact (resp. countably compact, resp.
Lindel\"{o}f) fibers. Once the concept of closed maps is extended to
convergence spaces, all the other notions extend as well in the obvious way.
As observed in \cite[section 10]{D.comp}, preservation of closed sets by a
map $f:(X,\xi )\rightarrow (Y,\tau )$ is equivalent to $\mathbb{F}_{1}$%
-compactness of the inverse map $f^{-}$ when $(X,\xi )$ is topological, but
not if $\xi $ is a general convergence. More precisely, calling a map $%
f:(X,\xi )\rightarrow (Y,\tau )$ \emph{adherent }\cite{D.comp} if 
\begin{equation*}
y\in \mathrm{adh}_{\tau }f(H)\Longrightarrow \mathrm{adh}%
_{\xi }H\cap f^{-}y\neq \emptyset ,
\end{equation*}%
we have:

\begin{lemma}
\begin{enumerate}
\item \label{lem:closed} A map $f:(X,\xi )\rightarrow (Y,\tau )$ is adherent
if and only if $f^{-}:(Y,\tau )\rightrightarrows (X,\xi )$ is an $\mathbb{F}%
_{1}$-compact relation;

\item If $f:(X,\xi )\rightarrow (Y,\tau )$ is adherent, then it is closed;

\item If $f:(X,\xi )\rightarrow (Y,\tau )$ is closed and if adherence of
sets are closed in $\xi $ (in particular if $\xi $ is a topology), then $f$
is adherent.
\end{enumerate}
\end{lemma}

\begin{proof}
$\left( 1\right) $ follows from the definition and is observed in \cite[%
section 10]{D.comp}.

(2) If $f(H)$ is not $\tau $-closed, then there exists $y\in \mathrm{adh}%
_{\tau }f(H)\diagdown f(H).$ Since $f$ is adherent, there exists $%
x\in \mathrm{adh}_{\xi }H\cap f^{-}y.$ But $x\notin H$ because $%
f(x)=y\notin f(H).$Therefore $H$ is not $\xi $-closed.

(3) is proved in \cite[Proposition 10.2]{D.comp} even if this proposition is
stated with a stronger assumption.
\end{proof}

Hence, a map $f:(X,\xi )\rightarrow (Y,\tau )$ with a domain in which
adherence of subsets are closed (in particular, a map with a topological
domain) is adherent if and only if it is closed if and only if $%
f^{-}:(Y,\tau )\rightrightarrows (X,\xi )$ is an $\mathbb{F}_{1}$-compact
relation. If the domain and range of a map are topological spaces, it is
well known that closedness of the map amounts to upper semicontinuity of the
inverse relation. It was observed (for instance in \cite{active}) that a
(multivalued) map is upper semicontinuous (u.s.c.) if and only if it is an $%
\mathbb{F}_{1}$-compact relation.

A surjection $f:X\rightarrow Y$ is $\mathbb{D}$\emph{-perfect }if it is
adherent with $\mathbb{D}$-compact fibers. In view of Corollary \ref%
{cor:fibers}, compact valued u.s.c. maps between topological spaces, known
as \emph{usco maps}, are compact relations. Another direct consequence of
Lemma \ref{lem:closed} and of Corollary \ref{cor:fibers} is:

\begin{theorem}
\label{th:Dperfect}Let $f:(X,\xi )\rightarrow (Y,\tau )$ be a surjection,
let $\mathbb{D}$ be an $\mathbb{F}_{1}$-composable class of filters, and let 
$\xi $ be a topology such that $\mathrm{adh}_{\xi }^{\natural }%
\mathbb{D\subset D}$. The following are equivalent:

\begin{enumerate}
\item $f$ is $\mathbb{D}$-perfect;

\item $f^{-}:Y\rightrightarrows X$ is $\mathbb{D}$-compact;

\item $f^{-}:Y\rightrightarrows X$ is $\left( \int \mathbb{D}\right) $%
-compact;

\item $f$ is $\left( \int \mathbb{D}\right) $-perfect.
\end{enumerate}
\end{theorem}

The equivalence between the first two points was first observed in \cite[%
Proposition 10.2]{D.comp} but erroneously stated for general convergences as
domain and range. Indeed, if $f:(X,\xi )\rightarrow (Y,\tau )$ is a
surjective map between two convergence spaces and if $f^{-}:\left( Y,\tau
\right) \rightrightarrows (X,\xi )$ is $\mathbb{D}$-compact, then $f$ is
adherent and has $\mathbb{D}$-compact fibers; if on the other hand $f$ is
adherent and has ($\mathbb{D}/\mathbb{F}_{1})$-compact fibers then $%
f^{-}:\left( Y,\tau \right) \rightrightarrows (X,\xi )$ is $\mathbb{D}$%
-compact. Hence, the two concepts are equivalent only when $\mathbb{D}$%
-compact sets are ($\mathbb{D}/\mathbb{F}_{1})$-compact in $\xi $, for
instance if $\mathrm{adh}_{\xi }^{\natural }\left( \mathbb{D}%
\right) \mathbb{\subset D}$ and $\xi $ is a $P$-diagonal convergence (in
particular if $\xi $ is a topology).

S. Dolecki offered in \cite{quest2} a unified treatment of various classes
of quotient maps and preservation theorems under such maps in the general
context of convergences. He extended the usual notions of quotient maps to
convergence spaces in the following way: a surjection $f:(X,\xi )\rightarrow
(Y,\tau )$ is $\mathbb{D}$\emph{-quotient} if 
\begin{equation}
y\in \mathrm{adh}_{\tau }\mathcal{H}\Longrightarrow f^{-}(y)\cap 
\mathrm{adh}_{\xi }f^{-}\mathcal{H\neq \emptyset },
\label{eq:Dquotient}
\end{equation}%
for every $\mathcal{H}\in \mathbb{D}(Y).$ When $\mathbb{D}$ is the class of
all (resp. countably based, principal, principal of closed sets) filters,
then continuous $\mathbb{D}$-quotient maps betweeen topological spaces are
exactly biquotient (resp. countably biquotient, hereditarily quotient,
quotient) maps. Now, I present a new characterization of $\mathbb{D}$%
-quotient maps as $\mathbb{D}$-compact relations, in this general context of
convergence spaces. As mentioned before, the category of convergence spaces
and continuous maps is topological, hence if $f:(X,\xi )\rightarrow Y,$
there exists the finest convergence ---called\emph{\ final convergence} and
denoted $f\xi $ --- on $Y$ making $f$ continuous. Analogously, if $%
f:X\rightarrow (Y,\tau ),$ there exists the coarsest convergence --- called 
\emph{initial convergence} and denoted $f^{-}\tau $ --- on $X$ making $f$
continuous. If $\tau $ is topological, so is $f^{-}\tau .$ In contrast, $%
f\xi $ can be non topological even when $\xi $ is topological.

\begin{theorem}
\label{th:Dquotient} Let $\mathbb{D}$ be an $\mathbb{F}_{1}$-composable
class of filters. Let $f:(X,\xi )\rightarrow (Y,\tau )$ be a surjection. The
following are equivalent:

\begin{enumerate}
\item $f:(X,\xi )\rightarrow (Y,\tau )$ is $\mathbb{D}$-quotient;

\item $\tau \geq \mathrm{Adh}_{\mathbb{D}}f\xi ;$

\item $f:\left( X,f^{-}\tau \right) \rightarrow (Y,f\xi )$ is a $\mathbb{D}$%
-compact relation.
\end{enumerate}
\end{theorem}

\begin{proof}
The equivalence $\left( 1\Longleftrightarrow 2\right) $ is \cite[Theorem 1.2]%
{quest2}.

$\left( 1\Longleftrightarrow 3\right) .$ Assume $f$ is $\mathbb{D}$-quotient
and let $x\in \lm_{f^{-}\tau }\mathcal{F}.$ Then $f(x)\in \lm_{\tau }f(%
\mathcal{F}),$ so that $f(x)\in \mathrm{adh}_{\tau }\mathcal{D}$ whenever $%
\mathcal{D\in }\mathbb{D}(Y)$ and $\mathcal{D}\#f(\mathcal{F}).$ By $\left( %
\ref{eq:Dquotient}\right) ,$ $f^{-}(f(x))\cap \mathrm{adh}_{\xi
}f^{-}\mathcal{D\neq \emptyset }$ so that $f(x)\in f\left( \mathrm{adh}%
_{\xi }f^{-}\mathcal{D}\right) .$ In view of \cite[Lemma 2.1]{D.comp}, $%
f(x)\in \mathrm{adh}_{f\xi }\mathcal{D}.$

Conversely, assume that $f:\left( X,f^{-}\tau \right) \rightarrow (Y,f\xi )$
is $\mathbb{D}$-compact and let $y\in \mathrm{adh}_{\tau }\mathcal{D}.$
There exists $\mathcal{G}\#\mathcal{D}$ such that $y\in \lm_{\tau }\mathcal{%
G}.$ By definition of $f^{-}\tau ,$ the filter $f^{-}\mathcal{G}$ is
converging to every point of $f^{-}y$ for $f^{-}\tau .$ In view of
Proposition \ref{pro:pointsuffice}, $ff^{-}\mathcal{G}$ is $\mathbb{D}$%
-compact at $\{y\}$ for $f\xi .$ Since $f$ is surjective, $ff^{-}\mathcal{G}%
\approx \mathcal{G}$ and $\mathcal{G}\#\mathcal{D}$ so that $y\in \mathrm{%
adh}_{f\xi }\mathcal{D}=f\left( \mathrm{adh}_{\xi }f^{-}\mathcal{D}\right)
, $ by \cite[Lemma 2.1]{D.comp}. Therefore, $f^{-}(y)\cap \mathrm{adh}%
_{\xi }f^{-}\mathcal{D\neq \emptyset }.$
\end{proof}

Notice that even if the map has topological range and domain, you need to
extend the notions to convergence spaces to obtain such a characterization.

\section{Compactly meshable filters and relations}

In view of the characterizations above of various types of maps as $\mathbb{D%
}$-compact relations, results of stability of $\mathbb{D}$-compactness of
filters under product would particularize to product theorems for $\mathbb{D}
$-compact spaces, but also for various types of quotient maps, for variants
of perfect and closed maps, for usc and usco maps. Product theorems for $%
\mathbb{D}$-compact filters and their applications is the purpose of the
companion paper \cite{myn.applofcompact}. A (complicated but extremely
useful) notion fundamental to this study of products is the following:

A filter $\mathcal{F}$ is $\mathbb{M}$-\emph{compactly} $\mathbb{J}$ \emph{to%
} $\mathbb{D}$ \emph{meshable }at $A,$ or $\mathcal{F}$ is an $\mathbb{M}$%
-compactly $\left( \mathbb{J}/\mathbb{D}\right) _{\#}$-filter at $A$, if 
\begin{equation*}
\mathcal{J}\in \mathbb{J},\mathcal{J}\#\mathcal{F}\Longrightarrow \exists 
\mathcal{D}\in \mathbb{D},\mathcal{D}\#\mathcal{J}\text{ and }\mathcal{D}%
\text{ is }\mathbb{M}\text{-compact at }A.
\end{equation*}

While the importance of this concept will be best highlighted by how it is
used in the companion paper \cite{myn.applofcompact}, I show here that the
notion of an $\mathbb{M}$-compactly $\left( \mathbb{J}/\mathbb{D}\right)
_{\#}$-filter is instrumental in characterizing a large number of classical
concepts.

The notion of total countable compactness was first introduced by Z. Frol%
\'{\i}k \cite{frolikpseudo} for a study of product of countably compact and
pseudocompact spaces and rediscovered under various names by several authors
(see \cite[p. 212]{vaughan}). A topological space $X$ is \emph{totally} 
\emph{countably compact} if every countably based filter has a finer
(equivalently, meshes a) compact countably based filter. The name comes from 
\emph{total }nets of Pettis\emph{. }Obviously, a topological space is
totally countably compact if and only if $\{X\}$ is compactly $\mathbb{F}%
_{\omega }$ to $\mathbb{F}_{\omega }$ meshable. In \cite{vaughan}, J.
Vaughan studied more generally under which condition a product of $\mathbb{D}
$-compact spaces is $\mathbb{D}$-compact, under mild conditions on the class
of filters $\mathbb{D}$. He used in particular the concept of a totally $%
\mathbb{D}$-compact space $X$, which amounts to $\{X\}$ being a compactly $%
\mathbb{(D}/\mathbb{D)}_{\#}$-filter.

On the other hand, Theorem \ref{th:accessible} can be completed by the
following immediate rephrasing of the notion of $\mathbb{M}$-compactly $(%
\mathbb{J}/\mathbb{D})_{\#}$-filters relative to a singleton in convergence
theoretic terms.

\begin{proposition}
\label{pro:local M-compactoidly}Let $\mathbb{D}$, $\mathbb{J}$ and $\mathbb{M%
}$ be three classes of filters, and let $\xi $ and $\theta $ be two
convergences on $X$. The following are equivalent:

\begin{enumerate}
\item $\theta \geq \mathrm{Adh}_{\mathbb{J}}\mathrm{Base}_{\mathbb{D}}%
\mathrm{Adh}_{\mathbb{M}}\xi ;$

\item $\mathcal{F}$ is an $\mathbb{M}$-compactly $(\mathbb{J}/\mathbb{D}%
)_{\#}$-filter at $\{x\}$ in $\xi $ whenever $x\in \lm_{\theta }\mathcal{F}%
. $

In particular, $\xi =\mathrm{Adh}_{\mathbb{M}}\xi $ is $(\mathbb{J}/\mathbb{%
D})$-accessible if and only if $\mathcal{F}$ is an $\mathbb{M}$-compactly $(%
\mathbb{J}/\mathbb{D})_{\#}$-filter at $\{x\}$ whenever $x\in \lm \mathcal{F%
}.$
\end{enumerate}
\end{proposition}

In view of Table 1, this applies to a variety of classical local topological
properties.

A \emph{relation} $R:(X,\xi )\rightrightarrows (Y,\tau )$ \emph{is} $\mathbb{%
M}$-\emph{compactly} $\left( \mathbb{J}/\mathbb{D}\right) $\emph{-meshable}
if 
\begin{equation*}
\mathcal{F}\underset{\xi }{\rightarrow }x\Longrightarrow R(\mathcal{F})\text{
is }\mathbb{M}\text{-compactly }\left( \mathbb{J}/\mathbb{D}\right) \text{%
-meshable at }Rx\text{ in }\tau .
\end{equation*}

\begin{theorem}
\label{th:Mquot+range}Let $\mathbb{M\subset J}$, let $\tau =\mathrm{Adh}%
_{\mathbb{M}}\tau $ and let $f:(X,\xi )\rightarrow (Y,\tau )$ be a
continuous surjection. The map $f$ is $\mathbb{M}$-quotient with $\mathbb{(J}%
/\mathbb{D)}$-accessible range if and only if $f:(X,f^{-}\tau )\rightarrow
(Y,f\xi )$ is an $\mathbb{M}$-compactly $\left( \mathbb{J}/\mathbb{D}\right) 
$-meshable relation.
\end{theorem}

\begin{proof}
Assume that $f$ is $\mathbb{M}$-quotient with $\mathbb{(J}/\mathbb{D)}$%
-accessible range and let $x\in \lm_{f^{-}\tau }\mathcal{F}.$ Then $%
y=f(x)\in \lm_{\tau }f(\mathcal{F}).$ Let $\mathcal{J}$ be a $\mathbb{J}$%
-filter such that $\mathcal{J}\#f(\mathcal{F}).$ Since $y\in \mathrm{adh}%
_{\tau }\mathcal{J}$ and $\tau $ is $\mathbb{(J}/\mathbb{D)}$%
-accessible, there exists a $\mathbb{D}$-filter $\mathcal{D}\#\mathcal{J}$
such that $y\in \lm_{\tau }\mathcal{D}.$ To show that $f(\mathcal{F})$ is $%
\mathbb{M}$-compactly $\left( \mathbb{J}/\mathbb{D}\right) $-meshable at $y$
in $f\xi ,$ it remains to show that $\mathcal{D}$ is $\mathbb{M}$-compact at 
$\{y\}$ for $f\xi ,$ that is, that $y\in \lm_{\mathrm{Adh}_{%
\mathbb{M}}f\xi }\mathcal{D},$ which follows from the $\mathbb{M}$%
-quotientness of $f.$

Conversely, assume that $f:(X,f^{-}\tau )\rightarrow (Y,f\xi )$ is an $%
\mathbb{M}$-compactly $\left( \mathbb{J}/\mathbb{D}\right) $-meshable
relation, and let $y\in \lm_{\tau }\mathcal{G}.$ Then $f^{-}(\mathcal{G})$
converges to any point $x\in f^{-}y$ for $f^{-}\tau .$ Therefore, $f(f^{-}%
\mathcal{G})$ is $\mathbb{M}$-compactly $\left( \mathbb{J}/\mathbb{D}\right) 
$-meshable at $\{y\}$ in $f\xi .$ Because $f$ is a surjection, $f(f^{-}%
\mathcal{G})=\mathcal{G}.$ Consider $\mathcal{M}\mathbb{\in M\subset J}$
such that $\mathcal{M}\#\mathcal{G}.$ There exists a $\mathbb{D}$-filter $%
\mathcal{D}\#\mathcal{M}$ which is $\mathbb{M}$-compact at $\{y\}$ in $f\xi
. $ Hence, $y\in \mathrm{adh}_{f\xi }\mathcal{M},$ so that $y\in
\lm_{\mathrm{Adh}_{\mathbb{M}}f\xi }\mathcal{G}.$ Therefore, $f$ is $%
\mathbb{M}$-quotient. Moreover, if $y\in \mathrm{adh}_{\tau }%
\mathcal{J}$ for a $\mathbb{J}$-filter $\mathcal{J},$ then there exists $%
\mathcal{G}\#\mathcal{J}$ such that $y\in \lm_{\tau }\mathcal{G}.$ By the
previous argument, $\mathcal{G}$ is $\mathbb{M}$-compactly $\left( \mathbb{J}%
/\mathbb{D}\right) $-meshable at $\{y\}$ in $f\xi .$ In particular, there
exists a $\mathbb{D}$-filter $\mathcal{D}\#\mathcal{J}$ which is $\mathbb{M}$%
-compact at $\{y\}$ in $f\xi .$ In other words, $y\in \lm_{\mathrm{Adh}%
_{\mathbb{M}}f\xi }\mathcal{D}.$ Since $f:(X,\xi )\rightarrow
(Y,\tau )$ is continuous, $\tau \leq f\xi $ so that $\mathrm{Adh}_{%
\mathbb{M}}\tau =\tau \leq \mathrm{Adh}_{\mathbb{M}}f\xi .$ Hence $%
y\in \lm_{\tau }\mathcal{D}$ and $\tau $ is $\mathbb{(J}/\mathbb{D)}$%
-accessible.
\end{proof}

\begin{center}
\begin{tabular}{|c|c|c|c|}
\hline
$\mathbb{M}$ & $\mathbb{J}$ & $\mathbb{D}$ & map $f$ as in Theorem \ref%
{th:Mquot+range} \\ \hline
$\mathbb{F}_{1}{}$ & $\mathbb{F}{}$ & $\mathbb{F}_{1}{}$ & hereditarily
quotient with finitely generated range \\ \hline
$\mathbb{F}_{1}{}$ & $\mathbb{F}_{1}{}$ & $\mathbb{F}_{\omega }{}$ & 
hereditarily quotient with Fr\'{e}chet range \\ \hline
$\mathbb{F}_{1}{}$ & $\mathbb{F}_{\omega }{}$ & $\mathbb{F}_{\omega }{}$ & 
hereditarily quotient with strongly Fr\'{e}chet range \\ \hline
$\mathbb{F}_{1}{}$ & $\mathbb{F}{}$ & $\mathbb{F}_{\omega }{}$ & 
hereditarily quotient with bisequential range \\ \hline
$\mathbb{F}_{1}$ & $\mathbb{F}$ & $\mathbb{F}$ & hereditarily quotient \\ 
\hline
$\mathbb{F}_{\omega }$ & $\mathbb{F}_{\omega }{}$ & $\mathbb{F}_{1}$ & 
countably biquotient with finitely generated range \\ \hline
$\mathbb{F}_{\omega }$ & $\mathbb{F}_{\omega }$ & $\mathbb{F}_{\omega }{}$ & 
countably biquotient with strongly Fr\'{e}chet range \\ \hline
$\mathbb{F}_{\omega }{}$ & $\mathbb{F}{}$ & $\mathbb{F}_{\omega }{}$ & 
countably biquotient with bisequential range \\ \hline
$\mathbb{F}_{\omega }{}$ & $\mathbb{F}{}$ & $\mathbb{F}{}$ & countably
biquotient \\ \hline
$\mathbb{F}{}$ & $\mathbb{F}{}$ & $\mathbb{F}_{1}{}$ & biquotient with
finitely generated range \\ \hline
$\mathbb{F}{}$ & $\mathbb{F}{}$ & $\mathbb{F}_{\omega }{}$ & biquotient with
bisequential range \\ \hline
$\mathbb{F}{}$ & $\mathbb{F}{}$ & $\mathbb{F}{}$ & biquotient \\ \hline
\end{tabular}
\end{center}

\begin{theorem}
\label{th:Mperfect+range}Let $\mathbb{M\subset J}$ and $\mathbb{D}$ be three
classes of filters, where $\mathbb{J}$ and $\mathbb{D}$ are $\mathbb{F}_{1}$%
-composable. Let $\tau =\mathrm{Adh}_{\mathbb{M}}\tau $ and let $%
\xi $ be a $P$-diagonal convergence such that $\mathrm{adh}_{\xi
}^{\natural }(\mathbb{M)\subset M}$ . Let $f:(X,\xi )\rightarrow (Y,\tau )$
be a continuous surjection. The map $f$ is $\mathbb{M}$-perfect with $%
\mathbb{(J}/\mathbb{D)}$-accessible range if and only if $f^{-}:(Y,\tau
)\rightrightarrows (X,\xi )$ is an $\mathbb{M}$-compactly $\left( \mathbb{J}/%
\mathbb{D}\right) $-meshable relation.
\end{theorem}

\begin{proof}
Assume that $f$ is $\mathbb{M}$-perfect with $\mathbb{(J}/\mathbb{D)}$%
-accessible range and let $y\in \lm_{\tau }\mathcal{G}.$ Consider a $%
\mathbb{J}$-filter $\mathcal{J}\#f^{-}\mathcal{G}.$ By $\mathbb{F}_{1}$%
-composability, $f(\mathcal{J})$ is a $\mathbb{J}$-filter. Moreover $f(%
\mathcal{J)}\#\mathcal{G}$ so that $y\in \mathrm{adh}_{\tau }f(%
\mathcal{J}).$ Since $\tau $ is $\mathbb{(J}/\mathbb{D)}$-accessible, there
exists a $\mathbb{D}$-filter $\mathcal{D}\#f(\mathcal{J})$ such that $y\in
\lm_{\tau }\mathcal{D}.$ In view of Corollary \ref{cor:fibers}, $%
f^{-}:(Y,\tau )\rightrightarrows (X,\xi )$ is $\mathbb{M}$-compact because $%
f $ is $\mathbb{M}$-perfect. Therefore, $f^{-}\mathcal{D}$ is $\mathbb{M}$%
-compact at $f^{-}y$ in $\xi .$ Moreover, $f^{-}\mathcal{D}\in \mathbb{D}$
because $\mathbb{D}$ is $\mathbb{F}_{1}$-composable and $f^{-}\mathcal{D}\#%
\mathcal{J}.$ Hence, $f^{-}\mathcal{G}$ is $\mathbb{M}$-compactly $\left( 
\mathbb{J}/\mathbb{D}\right) $-meshable at $f^{-}y$ in $\xi .$

Conversely, assume that $f^{-}:(Y,\tau )\rightrightarrows (X,\xi )$ is an $%
\mathbb{M}$-compactly $\left( \mathbb{J}/\mathbb{D}\right) $-meshable
relation. It is in particular an $\mathbb{M}$-compact relation because $%
\mathbb{M\subset J}.$ In view of Corollary \ref{cor:fibers}, $f$ is $\mathbb{%
M}$-perfect. Now assume that $y\mathrm{adh}_{\tau }\mathcal{J}$
where $\mathcal{J}\in \mathbb{J}.$ There exists $\mathcal{G}\#\mathcal{J}$
such that $y\in \lm_{\tau }\mathcal{G}.$ Therefore, $f^{-}\mathcal{G}$ is $%
\mathbb{M}$-compactly $\left( \mathbb{J}/\mathbb{D}\right) $-meshable at $%
f^{-}y$ in $\xi .$ The filter $f^{-}\mathcal{J}$ meshes with $f^{-}\mathcal{G%
}$ because $f$ is surjective, and is a $\mathbb{J}$-filter because $\mathbb{J%
}$ is $\mathbb{F}_{1}$-composable. Hence, there exists a $\mathbb{D}$-filter 
$\mathcal{D}\#f^{-}\mathcal{J}$ which is $\mathbb{M}$-compact at $f^{-}y$ in 
$\xi .$ By continuity of $f$, the filter $f(\mathcal{D})$ is $\mathbb{M}$%
-compact at $\{y\}$ in $\tau $ (Corollary \ref{cor:continuous}). In view of
Proposition \ref{pro:AdhD}, $y\in \lm_{\mathrm{Adh}_{\mathbb{M}%
}\tau }f(\mathcal{D}).$ Moreover, the filter $f(\mathcal{D})$ meshes with $%
\mathcal{J}$ and is a $\mathbb{D}$-filter, by $\mathbb{F}_{1}$-composability
of $\mathbb{D}.$ Since, $\tau =\mathrm{Adh}_{\mathbb{M}}\tau ,$ $%
y\in \lm_{\tau }f(\mathcal{D})$ and $\tau $ is $\mathbb{(J}/\mathbb{D)}$%
-accessible.
\end{proof}

\begin{center}
\begin{tabular}{|c|c|c|c|}
\hline
$\mathbb{M}$ & $\mathbb{J}$ & $\mathbb{D}$ & map $f$ as in Theorem \ref%
{th:Mperfect+range} \\ \hline
$\mathbb{F}_{1}{}$ & $\mathbb{F}{}$ & $\mathbb{F}_{1}{}$ & closed with
finitely generated range \\ \hline
$\mathbb{F}_{1}{}$ & $\mathbb{F}_{1}{}$ & $\mathbb{F}_{\omega }{}$ & closed
with Fr\'{e}chet range \\ \hline
$\mathbb{F}_{1}{}$ & $\mathbb{F}_{\omega }{}$ & $\mathbb{F}_{\omega }{}$ & 
closed with strongly Fr\'{e}chet range \\ \hline
$\mathbb{F}_{1}{}$ & $\mathbb{F}{}$ & $\mathbb{F}_{\omega }{}$ & closed with
bisequential range \\ \hline
$\mathbb{F}_{1}$ & $\mathbb{F}$ & $\mathbb{F}$ & closed \\ \hline
$\mathbb{F}_{\omega }$ & $\mathbb{F}_{\omega }{}$ & $\mathbb{F}_{1}$ & 
countably perfect with finitely generated range \\ \hline
$\mathbb{F}_{\omega }$ & $\mathbb{F}_{\omega }$ & $\mathbb{F}_{\omega }{}$ & 
countably perfect with strongly Fr\'{e}chet range \\ \hline
$\mathbb{F}_{\omega }{}$ & $\mathbb{F}{}$ & $\mathbb{F}_{\omega }{}$ & 
countably perfect with bisequential range \\ \hline
$\mathbb{F}_{\omega }{}$ & $\mathbb{F}{}$ & $\mathbb{F}{}$ & countably
perfect \\ \hline
$\mathbb{F}{}$ & $\mathbb{F}{}$ & $\mathbb{F}_{1}{}$ & perfect with finitely
generated range \\ \hline
$\mathbb{F}{}$ & $\mathbb{F}{}$ & $\mathbb{F}_{\omega }{}$ & perfect with
bisequential range \\ \hline
$\mathbb{F}{}$ & $\mathbb{F}{}$ & $\mathbb{F}{}$ & perfect \\ \hline
\end{tabular}
\end{center}

\bibliographystyle{amsplain}

\end{document}